\documentclass{amsart}
 
\newtheorem{maintheorem}{Theorem}
 
\newtheorem{maincorollary}[maintheorem]{Corollary}
\newtheorem{theorem}{Theorem}[section]
\newtheorem{lemma}[theorem]{Lemma}
\newtheorem{proposition}[theorem]{Proposition}
\newtheorem{corollary}[theorem]{Corollary}
\theoremstyle{definition}
\newtheorem{definition}[theorem]{Definition}

\newtheorem{remark}[theorem]{Remark}

\newtheorem{example}[theorem]{Example}
\numberwithin{equation}{theorem}
\newcommand{\SI}[1]{\textbf{SI(#1)}}
\newcommand{\MI}[1]{\textbf{MI(#1)}}
\newcommand{\boplus}{\textstyle \bigoplus}
\newcommand{\botimes}{\textstyle \bigotimes}

\newcommand{\HA}{H\!A}
\DeclareMathOperator{\adeg}{adeg}
\DeclareMathOperator{\Ext}{Ext}
\DeclareMathOperator{\End}{End}

\DeclareMathOperator{\Hom}{Hom}
\DeclareMathOperator{\im}{im}
\DeclareMathOperator{\Mer}{Mer}
\newcommand{\fm}{\mathfrak{m}}
\newcommand{\Z}{\mathbb{Z}}
\newcommand{\free}[1]{T\langle #1 \rangle}
\newcommand{\massey}[1]{\langle #1 \rangle}
\begin{document}

\title{$A$-infinity structure on Ext-algebras}

\author{D.-M. Lu, J. H. Palmieri, Q.-S. Wu and J. J. Zhang}

\address{(Lu) Department of Mathematics, Zhejiang University,
Hangzhou 310027, China}

\email{dmlu@zju.edu.cn}

\address{(Palmieri) Department of Mathematics, Box 354350, University of 
Washington, Seattle, WA 98195, USA}

\email{palmieri@math.washington.edu}

\address{(Wu) Institute of Mathematics, Fudan University, 
Shanghai, 200433, China}

\email{qswu@fudan.edu.cn}

\address{(Zhang) Department of Mathematics, Box 354350, University of 
Washington, Seattle, WA 98195, USA}

\email{zhang@math.washington.edu}

\begin{abstract}
Let $A$ be a connected graded algebra and let $E$ denote its
Ext-algebra $\boplus_i \Ext^i_A(k_A,k_A)$. There is a natural
$A_\infty$-structure on $E$, and we prove that this structure is
mainly determined by the relations of $A$. In particular, the
coefficients of the $A_{\infty}$-products $m_n$ restricted to the
tensor powers of $\Ext^1_A(k_A,k_A)$ give the coefficients of the
relations of $A$.  We also relate the $m_n$'s to Massey products.
\end{abstract}



\maketitle

\section*{Introduction}
\label{sect0}

The notions of $A_\infty$-algebra and $A_\infty$-space were 
introduced by Stasheff in the 1960s \cite{St}. Since then,
more and more theories involving $A_\infty$-structures (and its
cousins, $E_\infty$, $L_\infty$, and $B_\infty$) have 
been discovered in several areas of mathematics and physics. 
Kontsevich's talk \cite{Ko} at the ICM 1994 on categorical 
mirror symmetry has had an influence in developing 
this subject. The use of $A_\infty$-algebras in noncommutative 
algebra and the representation theory of algebras was introduced 
by Keller \cite{Ke1, Ke2, Ke3}. Recently the authors of this paper used 
the $A_\infty$-structure on the Ext-algebra $\Ext^*_A(k,k)$ to study 
the non-Koszul Artin-Schelter regular algebras $A$ of global 
dimension four \cite{LP3}. The information about the higher 
multiplications on $\Ext^*_A(k,k)$ is essential and very effective 
for the work \cite{LP3}. 

Throughout let $k$ be a commutative base field. The definition of an
$A_{\infty}$-algebra will be given in Section 1. Roughly speaking, an
$A_\infty$-algebra is a graded vector space $E$ equipped with a
sequence of ``multiplications'' $(m_1, m_2, m_3, \cdots)$: $m_{1}$ is
a differential, $m_2$ is the usual product, and the higher $m_{n}$'s
are homotopies which measure the degree of associativity of $m_{2}$.
An associative algebra $E$ (concentrated in degree 0) is an
$A_\infty$-algebra with multiplications $m_n=0$ for all $n\neq 2$, so
sometimes we write an associative algebra as $(E, m_2)$.  A
differential graded (DG) algebra $(E,d)$ has multiplication $m_2$ and
derivation $m_1=d$; this makes it into an $A_\infty$-algebra with
$m_{n}=0$ for $n \geq 3$ and so it could be written as $(E, m_1,
m_2)$.

Let $A$ be a connected graded algebra, and let $k_A$ be the right
trivial $A$-module $A/A_{\geq 1}$. The Ext-algebra $\boplus_{i\geq
0}\Ext^i_A(k_A,k_A)$ of $A$ is the homology of a DG algebra, and hence
by a theorem of Kadeishvili, it is equipped with an $A_\infty$-algebra
structure. We use $\Ext^*_A(k_A,k_A)$ to denote both the usual
associative Ext-algebra and the Ext-algebra with its
$A_\infty$-structure. By \cite[Ex. 13.4]{LP1} there is a graded
algebra $A$ such that the associative algebra $\Ext^*_A(k_A,k_A)$ does
not contain enough information to recover the original algebra $A$; on
the other hand, the information from the $A_\infty$-algebra
$\Ext^*_A(k_A,k_A)$ is sufficient to recover $A$. This is the point of
Theorem \ref{thmA} below, and this process of recovering the algebra
from its Ext-algebra is one of the main tools used in \cite{LP3}.

We need some notation in order to state the theorem.
Let $\fm = \boplus_{i \geq 1} A_{i}$ be the augmentation ideal of $A$.
We say that a graded vector space $V=\boplus V_{i}$ is \emph{locally
finite} if each $V_{i}$ is finite-dimensional. We write the graded
$k$-linear dual of $V$ as $V^{\#}$.  As our notation has so far
indicated, we use subscripts to indicate the grading on $A$ and
related vector spaces.  For example,
\[
(\fm^{\otimes m})_{s} = \boplus_{i_{1}+\dotsb +i_{m}=s} A_{i_{1}}
\otimes \dotsb \otimes A_{i_{m}}.
\]
Also, the grading on $A$ induces a bigrading on $\Ext$.  We write the
usual, homological, grading with superscripts, and the second,
induced, grading with subscripts.

Let $Q=\fm /\fm^{2}$ be the graded vector space of generators of
$A$.  The relations in $A$ naturally sit inside the tensor algebra on
$Q$.  In Section 4 we choose a vector space embedding of each graded
piece $A_{s}$ into the tensor algebra on $Q$: a map
\[
A_{s} \hookrightarrow (\boplus_{m \geq 1} Q^{\otimes m})_{s},
\]
which splits the multiplication map, and this choice affects how we
choose the minimal generating set of relations.  See Lemma \ref{lem5.2}
and the surrounding discussion for more details.

\begin{maintheorem}\label{thmA}
Let $A$ be a connected graded locally finite algebra, and let $E$ be
the $A_\infty$-algebra $\Ext^*_A(k_A,k_A)$.  Let $Q=\fm /\fm^{2}$ be
the graded vector space of generators of $A$.  Let $R=\boplus_{s\geq
2} R_s$ be a minimal graded space of relations of $A$, with $R_{s}$
chosen so that $R_s \subset \boplus_{1 \leq i \leq s-1} Q_{i} \otimes
A_{s-i} \subset (\boplus_{m \geq 2} Q^{\otimes m})_{s}$. For each $n
\geq 2$ and $s \geq 2$, let $i_{s}: R_s\to (\boplus_{m \geq 2}
Q^{\otimes m})_{s}$ be the inclusion map and let $i_{s}^{n}$ be the
composite
\[
R_{s} \xrightarrow{i_{s}} (\boplus_{m \geq 2} Q^{\otimes m})_{s}
\rightarrow (Q^{\otimes n})_{s}.
\]
Then in any degree $-s$, the multiplication $m_n$ of $E$ restricted
to $(E^1)^{\otimes n}_{s}$ is equal to the map
\[
(i_{s}^{n})^{\#}: \left((E^1)^{\otimes n}\right)_{-s} =
\left((Q^{\otimes n})_{s}\right)^{\#} \longrightarrow R^{\#}_s \subset
E^2_{-s}.
\]
\end{maintheorem}

In plain English, the multiplication maps $m_{n}$ on classes in
$\Ext^{1}_{A}(k_{A},k_{A})$ are determined by the relations in the
algebra $A$.

Note that the space $Q$ of generators need not be finite-dimensional
-- it only has to be finite-dimensional in each grading.  Thus it
applies to infinitely generated algebras like the Steenrod algebra.

The authors originally announced the result in the following special
case; this was used heavily in \cite{LP3}.

\begin{maincorollary}
[Keller's higher-multiplication theorem in the connected graded case]
\label{corB}
Let $A$ be a graded algebra, finitely generated in degree 1, and let
$E$ be the $A_\infty$-algebra $\Ext^*_A(k_A,k_A)$.  Let
$R=\boplus_{n\geq 2} R_n$ be a minimal graded space of relations of
$A$, chosen so that $R_n\subset A_1\otimes A_{n-1} \subset
A_1^{\otimes n} $. For each $n \geq 2$, let $i_{n}: R_n\to
A_1^{\otimes n}$ be the inclusion map and let $i_{n}^{\#}$ be its
$k$-linear dual. Then the multiplication $m_n$ of $E$ restricted to
$(E^1)^{\otimes n}$ is equal to the map
\[
i_{n}^{\#}: (E^1)^{\otimes n}=(A_1^{\#})^{\otimes n}\longrightarrow 
R^{\#}_n\subset E^2.
\]
\end{maincorollary}

Keller has the same result for a different class of algebras; indeed,
his result was the inspiration for Theorem \ref{thmA}.  His result
applies to algebras the form $k\Delta/I$ where $\Delta$ is a finite
quiver and $I$ is an admissible ideal of $k\Delta$; this was stated in
\cite[Proposition 2]{Ke4} without proof.  This class of algebras
includes those in Corollary \ref{corB}, but since the algebra $A$ in
Theorem \ref{thmA} need not be finitely generated, that theorem is not
a special case of Keller's result.  A version of Corollary \ref{corB}
was also proved in a recent paper by He and Lu \cite{HL} for ${\mathbb
N}$-graded algebras $A=A_0\oplus A_1\oplus \cdots$ with $A_0=k^{\oplus
n}$ for some $n\geq 1$, and which are finitely generated by $A_{0}
\oplus A_{1}$.  Their proof was based on the one here (see \cite[page
356]{HL}).

Here is an outline of the paper. We review the definitions of
$A_\infty$-algebras and Adams grading in Section \ref{sect1}.  In
Section \ref{sect2} we discuss Kadeishvili's and Merkulov's results
about the $A_{\infty}$-structure on the homology of a DG algebra.  In
Section \ref{sect3} we use Merkulov's construction to show that
the $A_{\infty}$-multiplication maps $m_{n}$ compute Massey products,
up to a sign -- see Theorem \ref{thm-massey} for details.  In Section
\ref{sect4} the bar construction is described: this is a DG algebra
whose homology is Ext, and so leads to an $A_{\infty}$-structure on
Ext algebras.  Then we give a proof of Theorem \ref{thmA} in Section
\ref{sect5}, and in Section \ref{sect6} we give a few examples.

This paper began as an appendix in \cite{LP3}. 

\section{Definitions} 
\label{sect1}

In this section we review the definition of an $A_{\infty}$-algebra
and discuss grading systems. Other basic material about 
$A_\infty$-algebras can be found in Keller's paper \cite{Ke3}.
Some examples of $A_\infty$-algebras related to ring theory were 
given in \cite{LP1}. Here is Stasheff's definition.

\begin{definition}
\label{defn1.1} \cite{St}
An {\it $A_\infty$-algebra} over a base field $k$ is a 
$\mathbb Z$-graded vector space
\[
A=\boplus_{p\in {\mathbb Z}} A^p
\]
endowed with a family of graded $k$-linear maps
\[
m_n: A^{\otimes n} \to A, \quad n\geq 1,
\]
of degree $2-n$ satisfying the following {\it Stasheff identities}:
\begin{equation}\label{SI}
\sum (-1)^{r+st} m_u(id^{\otimes r}\otimes m_s \otimes id^{\otimes
t})=0
\tag*{\SI{n}}
\end{equation}
for all $n\geq 1$, where the sum runs over all decompositions 
$n=r+s+t$ ($r,t\geq 0$ and $s\geq 1$), and where $u=r+1+t$. Here, 
$id$ denotes the identity map of $A$. Note that when these formulas 
are applied to elements, additional signs appear due to the Koszul 
sign rule. Some authors also use the terminology {\it strongly 
homotopy associative algebra} (or {\it sha algebra}) for 
$A_\infty$-algebra.
\end{definition}

The degree of $m_1$ is $1$ and the identity \SI{1} is $m_1 m_1=0$. This
says that $m_1$ is a differential of $A$. The identity \SI{2} is
\[
m_1 m_2= m_2(m_1\otimes id+id\otimes m_1)
\]
as maps $A^{\otimes 2} \to A$. So the differential
$m_1$ is a graded derivation with respect to $m_2$. Note that $m_2$
plays the role of multiplication although it may not be associative.
The degree of $m_2$ is zero. The identity \SI{3} is
\[
m_2(id\otimes m_2-m_2\otimes id)
=m_1m_3+ m_3(m_1\otimes id \otimes id+ id\otimes m_1\otimes
id + id\otimes id \otimes m_1)
\]
as maps $A^{\otimes 3}\to A$. If either $m_1$ or $m_3$ is zero, then
$m_2$ is associative. In general, $m_2$ is associative up to a
chain homotopy given by $m_3$.

When $n\geq 3$, the map $m_n$ is called a {\it higher
multiplication}. We write an $A_\infty$-algebra $A$ as $(A, m_1, m_2,
m_3, \cdots)$ to indicate the multiplications $m_i$.  We also assume
that every $A_\infty$-algebra in this paper contains an identity
element $1$ with respect to the multiplication $m_2$ that satisfies
the following {\it strictly unital condition}:

\bigskip

If $n\neq 2$ and $a_i=1$ for some $i$, then $m_n(a_1,\cdots, a_n)=0$. 

\bigskip
\noindent
In this case, $1$ is called the {\it strict unit} or {\it identity} of $A$.

We are mainly interested in graded algebras and their
Ext-algebras.  The grading appearing in a graded algebra may be
different from the grading appearing in the definition of the
$A_\infty$-algebra.  We introduce the Adams grading for an
$A_\infty$-algebra, as follows. Let $G$ be an abelian group.  (In this
paper, $G$ will always be free abelian of finite rank.)  Consider a
bigraded vector space
\[
A=\boplus_{p\in \Z, i\in G} A^p_i
\]
where the upper index $p$ is the grading appearing in Definition
\ref{defn1.1}, and the lower index $i$ is an extra grading, called
the \emph{$G$-Adams grading}, or {\it Adams grading} if $G$ is
understood.  We also write
\[
A^p=\boplus_{i\in G}A^p_i \quad\text{and}\quad A_i=
\boplus_{p\in \Z}A^p_i.
\]
The degree of a nonzero element in $A^p_i$ is $(p,i)$, and the second
degree is called the {\it Adams degree}.  For an $A_{\infty}$-algebra
$A$ to have an Adams grading, the map $m_n$ in Definition
\ref{defn1.1} must be of degree $(2-n,0)$: each $m_n$ must preserve
the Adams grading. When $A$ is an associative $G$-graded algebra
$A=\boplus_{i\in G} A_i$, we view $A$ as an $A_\infty$-algebra (or a
DG algebra) concentrated in degree 0, viewing the given grading on 
$A$ as the Adams grading.  The Ext-algebra of a graded algebra is 
bigraded; the grading inherited from the graded algebra is the 
Adams grading, and we keep using the lower index to denote the 
Adams degree.

Assume now that $G=\Z$, since we are mainly interested in this case.
Write 
\[
A^{\geq n}=\boplus_{p\geq n}A^p \quad\text{and}\quad A_{\geq n}=
\boplus_{i\geq n}A_i,
\]
and similarly for $A^{\leq n}$ and $A_{\leq n}$.
An
$A_\infty$-algebra $A$ with a $\Z$-Adams grading
is called {\it Adams connected} if (a) $A_0=k$, (b) $A=A_{\geq 0}$ 
or $A=A_{\leq 0}$, and (c) $A_i$ is finite-dimensional for all $i$. 
When $G=\Z\times G_0$, we define {\it Adams connected} 
in the same way after omitting the $G_0$-grading. If $A$ is a connected 
graded algebra which is finite-dimensional in each degree, then it is 
Adams connected when viewed as an $A_\infty$-algebra concentrated in 
degree 0. 

The following result is a consequence of Theorem \ref{thmA}, and it
will be proved at the end.  There might be several quasi-isomorphic
$A_{\infty}$-structures on $E:=\Ext_{A}^{*}(k_{A}, k_{A})$; we call
these different structures \emph{models} for the quasi-isomorphism
class of $E$.

\begin{proposition}\label{prop1.2} 
Let $A$ be a $\Z\oplus G$-Adams graded algebra, such that with respect
to the $\Z$-grading, $A$ is locally finitely generated.  Then there is
an $A_{\infty}$-model for $E$ such that the multiplications $m_n$ in
Theorem \ref{thmA} preserve the $\Z\oplus G$ grading.
\end{proposition}

\section{Kadeishvili's theorem and Merkulov's construction}
\label{sect2}

Let $A$ and $B$ be two $A_\infty$-algebras. A {\it morphism} of 
$A_\infty$-algebras $f: A\to B$ is a family of $k$-linear graded maps
\[
f_n: A^{\otimes n}\to B
\]
of degree $1-n$ satisfying the following 
{\it Stasheff morphism identities}:
\[
\sum (-1)^{r+st} f_u(id^{\otimes r}\otimes m_s\otimes 
id^{\otimes t})=\sum (-1)^{w} m_q(f_{i_1}\otimes f_{i_2}
\otimes \cdots \otimes f_{i_q})
\tag*{\MI{n}}
\]
for all $n\geq 1$, 
where the first sum runs over all decompositions $n=r+s+t$ with 
$s\geq 1$ and $r,t\geq 0$, where $u=r+1+t$, and the second sum runs 
over all $1\leq q\leq n$ and all decompositions $n=i_1+\cdots + i_{q}$
with all $i_s\geq 1$.  The sign on the right-hand side is given by
\[
w=(q-1)(i_1 -1) + (q-2) (i_2-1)+ \cdots + 2(i_{q-2}-1)+ (i_{q-1}-1).
\]

When $A$ and $B$ have a strict unit (as we always assume), 
an $A_\infty$-morphism is also required to satisfy the following extra 
{\it unital morphism conditions}:
\[
f_1(1_A)=1_B
\]
where $1_A$ and $1_B$ are the strict units of $A$ and $B$ respectively, and
\[
f_n(a_1 \otimes \cdots \otimes a_n)= 0 
\]
if $n\geq 2$ and $a_i=1_A$ for some $i$.

If $A$ and $B$ have Adams gradings indexed by the same group, then 
the maps $f_i$ are required to preserve the Adams degree.

A morphism $f$ is called a {\it quasi-isomorphism} if $f_1$ is a 
quasi-isomorphism.  A morphism is {\it strict} if $f_i=0$ for all 
$i\neq 1$. The \emph{identity morphism} is the strict morphism $f$ 
such that $f_1$ is the identity of $A$. When $f$ is a strict 
morphism from $A$ to $B$, then the identity \MI{n} becomes
\[
f_1 m_n=m_n(f_1\otimes \cdots \otimes f_1).
\]
A morphism $f=(f_i)$ is called a {\it strict isomorphism} if it is
strict with $f_1$ a vector space isomorphism.

Let $A$ be an $A_\infty$-algebra.  Its cohomology ring is defined 
to be
\[
\HA:=\ker m_1/\im m_1.
\] 
The following result, due to Kadeishvili \cite{Ka}, is a basic and
important property of $A_\infty$-algebras.

\begin{theorem}
\label{thm2.1} \cite{Ka}
Let $A$ be an $A_\infty$-algebra and let $\HA$ be the cohomology ring 
of $A$. There is an $A_\infty$-algebra structure on $\HA$ with $m_1=0$
and $m_{2}$ induced by the multiplication on $A$,
constructed from the $A_\infty$-structure of $A$, such that there is 
a quasi-isomorphism of $A_\infty$-algebras $\HA\to A$ lifting the 
identity of $\HA$.  This $A_\infty$-algebra structure on $\HA$ is unique 
up to quasi-isomorphism.
\end{theorem}

Kadeishvili's construction is very general. We would like to describe
some specific $A_\infty$-structures that we can work with.  Merkulov
constructed a special class of higher multiplications for $\HA$ in
\cite{Me1}, in which the higher multiplications can be defined
inductively; this way, the $A_\infty$-structure can be described
more explicitly, and hence used more effectively.  For our purposes we
will describe a special case of Merkulov's construction, assuming that
$A$ is a DG algebra.

Let $A$ be a DG algebra with differential $\partial$ and 
multiplication $\cdot$. Denote by $B^n$ and $Z^n$
the coboundaries and cocycles of $A^n$, respectively. Then there
are subspaces $H^n$ and $L^n$ such that 
\[
Z^n=B^n\oplus H^n
\]
and 
\begin{equation}
\label{2.1.1}
A^n=Z^n\oplus L^n=B^n\oplus H^n\oplus L^n.
\end{equation}
We will identify $\HA$ with $\boplus_{n} H^n$, or embed $\HA$ into $A$ 
by cocycle-sections $H^n\subset A^n$. {\it There are many different 
choices of $H^n$ and $L^n$.} 

Note that if $A$ has an Adams grading, then the decompositions above
will be chosen to respect the Adams grading, and all
maps constructed below will preserve the Adams grading.

Let $p=Pr_H: A\to A$ be a projection to $H:=\boplus_n H^n$, and 
let $G: A\to A$ be a homotopy from $id_A$ to $p$. Hence 
we have $id_A-p=\partial G+G\partial$.  The map $G$ is not unique, and
we want to choose $G$ carefully, so we define it as follows:
for every $n$, $G^n: A^n\to A^{n-1}$ is the map which satisfies
\begin{itemize}
\item 
$G^n=0$ when restricted to $L^n$ and $H^n$, and 
\item
$G^n=(\partial^{n-1}|_{L^{n-1}})^{-1}$ when restricted to $B^n$.
\end{itemize}
So the image of $G^n$ is $L^{n-1}$. 
It follows that $G^{n+1}\partial^n=Pr_{L^n}$ and 
$\partial^{n-1}G^n=Pr_{B^n}$.

Define a sequence of linear maps $\lambda_n: A^{\otimes n}\to A$ of
degree $2-n$ as follows.  There is no map $\lambda_{1}$, but we
formally define the ``composite'' $G \lambda_{1}$ by $G \lambda_{1} =
-id_{A}$.  $\lambda_2$ is the multiplication of $A$, namely, 
$\lambda_2(a_1\otimes a_2)=a_1\cdot a_2$.
For $n\geq 3$,
$\lambda_n$ is defined by the recursive formula
\begin{equation}
\label{2.1.2}
\lambda_n=\sum_{\substack{s+t=n, \\ s,t\geq 1}} (-1)^{s+1} 
\lambda_2[G\lambda_s\otimes G\lambda_t].
\end{equation}

We abuse notation slightly, and use $p$ to denote both the map $A
\rightarrow A$ and also (since the image of $p$ is $\HA$) the map
$A \rightarrow \HA$; we also use $\lambda_{i}$ both for the map
$A^{\otimes i} \rightarrow A$ and for its restriction $(\HA)^{\otimes
i} \rightarrow A$ to $\HA^{\otimes i}$.

Merkulov reproved Kadeishvili's result in \cite{Me1}.

\begin{theorem}
\label{thm2.2} \cite{Me1}
Let $m_i=p\lambda_i$. Then $(\HA, m_2,m_3, \cdots)$ is an 
$A_\infty$-algebra.
\end{theorem}

We can also display the quasi-isomorphism between $\HA$ and $A$ directly.

\begin{proposition}
\label{prop2.3} 
Let $\{\lambda_n\}$ be defined as above. For $i \geq 1$ let $f_i=-G 
\lambda_i: (\HA)^{\otimes i}\to A$, and for $i \geq 2$ let $m_i=p \lambda_i:
(\HA)^{\otimes i}\to \HA$. Then $(\HA, m_2, m_3, \cdots)$ is an 
$A_\infty$-algebra
and $f:=\{f_i\}$ is a quasi-isomorphism of $A_\infty$-algebras.
\end{proposition}

\begin{proof} 
This construction of $\{m_i\}$ and $\{f_i\}$ is a special case of
Kadeishvili's construction.
\end{proof}

Any $A_\infty$-algebra constructed as in Theorem \ref{thm2.2} and 
Proposition \ref{prop2.3} 
is called a {\it Merkulov model} of $A$, denoted by $H_{\Mer}A$. 
The particular model depends on the decomposition \eqref{2.1.1}, but 
all Merkulov models of $A$ are quasi-isomorphic to each other. 
If $A$ has an Adams grading, then by construction all maps $m_i$ and 
$f_i$ preserve the Adams degree. 

Next we consider the unital condition. 

\begin{lemma}
\label{lem2.4} Suppose $H^0$ is chosen to contain the unit element 
of $A$. Then $H_{\Mer}A$ satisfies the strictly unital condition, 
and the morphism $f=\{f_i\}$ satisfies the unital morphism conditions.
\end{lemma}

\begin{proof} First of all, $1\in H^0$ is a unit with respect 
to $m_2$.  We use induction on $n$ to show the following, for $n \geq 3$:

$(a)_n$: $\ $ $f_{n-1}(a_1 \otimes \dotsb \otimes a_{n-1})=0$ if $a_i=1$ 
for some $i$.

$(b)_n$: $\ $ $\lambda_n(a_1 \otimes \dotsb \otimes a_n)\in 
L:=\boplus_n L^n$ if $a_i=1$ for some $i$.

$(c)_n$: $\ $ $m_n(a_1 \otimes \dotsb \otimes a_n)=0$ if $a_i=1$ for 
some $i$.

\noindent
The strictly unital condition is $(c)_n$.  The unital morphism
condition is $(a)_n$.

We first prove $(a)_3$. For $a\in H$, 
\[ 
f_2(1 \otimes a)=-G\lambda_2(1 \otimes a)=-G(a)=0,
\]
since $G|_{H}=0$.  Similarly, $f_2(a \otimes 1)=0$.  This proves $(a)_3$. 
Now suppose for some $n \geq 3$ that $(a)_i$ holds for all 
$3 \leq i \leq n$.  By definition,
\[ 
\lambda_n=\sum_{s=1}^{n-1} (-1)^{s+1} \lambda_2( f_s\otimes f_{n-s}). 
\]
If $a_1=1$, $(a)_n$ implies that 
\[ 
\lambda_n(a_1 \otimes \dotsb \otimes a_n)=f_{n-1}(a_2 \otimes \dotsb
\otimes a_n)\in L. 
\]
Similarly, if $a_n=1$, we have $\lambda_n(a_1 \otimes \dotsb \otimes a_n)
\in L$. If $a_i=1$ for $1<i<n$, then $\lambda_n(a_1 \otimes \dotsb 
\otimes a_n)=0$. Therefore $(a)_i$ for $i \leq n$ implies $(b)_n$. 
Since $p(L)=0$, $(c)_n$ follows from $(b)_n$. Since $G(L)=0$, 
$(a)_{n+1}$ follows from $(b)_n$.  Induction completes the proof.
\end{proof}

\begin{lemma}
\label{lem2.5}
Let $(A,\partial)$ be a DG algebra and let $e\in A^0$ be an 
idempotent such that 
$\partial(e)=0$. Let $D=eAe$ and $C=(1-e)A+A(1-e)$.
\begin{enumerate}
\item If $HC=0$, then 
we can choose Merkulov models so that $H_{\Mer}A$ is strictly 
isomorphic to $H_{\Mer}D$. As a consequence
$A$ and $D$ are quasi-isomorphic as $A_\infty$-algebras.  
\item If moreover $\HA$ is Adams connected, then $H^0_{\Mer}A$ and $H^0_{\Mer}D$ 
in part (a) can be chosen to contain the unit element.
\end{enumerate}  
\end{lemma}

\begin{proof} 
First of all, $D$ is a sub-DG algebra of $A$ with identity $e$. 
Since $A=D\oplus C$ as chain complexes,
the group of coboundaries $B^{n}$ decomposes as $B^{n} = B^{n}_{D} \oplus
B^{n}_{C}$, where $B^{n}_{D} = B^{n} \cap D$ and $B^{n}_{C} = B^{n}
\cap C$.  Since $HC=0$, we can choose $H$ and $L$ so that they
decompose similarly (with $H_{C}=0$), giving the following direct sum
decompositions: 
\begin{gather*} 
A^n=D^n\oplus C^n= (B^n_D\oplus H^n_D\oplus L^n_D)\oplus 
(B^n_C\oplus L^n_C), \\
A^n=B^n\oplus H^n\oplus L^n=(B^n_D\oplus B^n_C)\oplus H^n_D\oplus
(L^n_D\oplus L^n_C).
\end{gather*}
It follows from the construction before Theorem \ref{thm2.2} that 
$H_{\Mer}A=H_{\Mer}D$. We choose $H^0_D$ to contain $e$. By Lemma 
\ref{lem2.4}, $e$ is the strict unit of $H_{\Mer}D$; hence $e$ is 
the strict unit of $H_{\Mer}A$, but note that the unit $1$ of $A$ 
may not be in $\HA$.  

Now suppose $\HA$ is Adams connected with unit $u$. Let $H^0=ku\oplus 
H^0_{\geq 1}$ (or $H^0=ku\oplus H^0_{\leq -1}$ if negatively connected 
graded). Replace $H^0$ by $k1\oplus H^0_{\geq 1}$ and keep the 
other subspaces $B^n$, $H^n$, and $L^n$ the same. Let $\overline{H_{\Mer}A}$ 
denote the new Merkulov model with the new choice of $H^0$. Then 
by Lemma \ref{lem2.4}, $1$ is the strict unit of $\overline{H_{\Mer}A}$. 
By construction, we have $(\overline{H_{\Mer}A})_{\geq 1}=
(H_{\Mer}A)_{\geq 1}$ as $A_\infty$-algebras without unit. By the 
unital condition, we see that $\overline{H_{\Mer}A}$ is strictly 
isomorphic to $H_{\Mer}A$. 
\end{proof}

\section{Massey products}
\label{sect3}

It is common to view $A_{\infty}$-algebras as algebras which are
strongly homotopy associative: not associative on the nose, but
associative up to all higher homotopies, as given by the $m_{n}$'s.
Any $A_{\infty}$-algebra in which the differential $m_{1}$ is zero,
such as the cohomology of a DG algebra, is strictly associative,
though, and in such a case, it is natural to wonder about the role of
the higher multiplications.  On the other hand, the cohomology of a DG
algebra is the natural setting for Massey products.  With Merkulov's
construction in hand, we give a proof of a folk theorem which connects
the higher multiplication maps with Massey products: we prove that
they are essentially the same, up to a sign.  We start by reviewing
Massey products.  We use the sign conventions from May \cite{Ma}; see
also Ravenel \cite[A1.4]{Ra}.

If $a$ is an element of a DG algebra $A$, we write $\overline{a}$ for
$(-1)^{1+\deg a}a$.  (This notation helps to keep some formulas simple.)

The length two Massey product $\massey{\alpha_{1}, \alpha_{2}}$ is the
ordinary product: $\massey{\alpha_{1}, \alpha_{2}} = \alpha_{1}
\alpha_{2}$.  (For consistency with the higher products, one could
also define $\massey{\alpha_{1}, \alpha_{2}}$ as being the set
$\{\alpha_{1}, \alpha_{2} \}$, but we do not take this point of
view.).

The Massey triple product is defined as follows: suppose given classes
$\alpha_{1}, \alpha_{2}, \alpha_{3} \in HA$ which are represented by
cocycles $a_{01}, a_{12}, a_{23} \in A$, respectively, and suppose
that $\alpha_{1} \alpha_{2} = 0 = \alpha_{2} \alpha_{3}$.  Then there
are cochains $a_{02}$ and $a_{13}$ so that
$\partial(a_{02})=\overline{a}_{01} a_{12}$ and
$\partial(a_{13})=\overline{a}_{12} a_{23}$.  Then
\[
\overline{a}_{02} a_{23} + \overline{a}_{01} a_{13}
\]
is a cocycle, and so represents a cohomology class.  One can choose
different cochains for $a_{02}$ and $a_{13}$: one can replace $a_{02}$
with $a_{02}+z$ for any cocycle $z$, for instance, and this can
produce a different cohomology class.  The \emph{length 3 Massey
product} $\massey{\alpha_{1}, \alpha_{2}, \alpha_{3}}$ is the set of
cohomology classes which arise from all such choices of $a_{02}$ and
$a_{13}$.

More generally, for any $n \geq 3$, the length $n$ Massey product is
defined as follows.  Suppose that we have cohomology classes
$\alpha_{i}$ for $1 \leq i \leq n$.  Suppose that whenever $i<j$ and $j-i<n-1$,
each length $j-i+1$ Massey product $\massey{\alpha_{i}, \dotsc,
\alpha_{j}}$ is defined and contains zero.  Then the \emph{length $n$
Massey product} $\massey{\alpha_{1}, \dotsc, \alpha_{n}}$ exists and
is defined as follows: using induction on $j-i+1$, one defines
cochains $a_{ij}$ as follows: $a_{i-1,i}$ is a cocycle representing
the cohomology class $\alpha_{i}$.  Given $a_{km}$ for all $k<m$ with
$m-k+1<j-i+1$, choose $a_{ij}$ so that
\[
\partial(a_{ij}) = \sum_{i<k<j} \overline{a}_{ik} a_{kj}.
\]
Then $\massey{\alpha_{1}, \dotsc, \alpha_{n}}$ is the set of cohomology
classes represented by cocycles of the form
\[
\sum_{0 < i < n} \overline{a}_{0i} a_{in}.
\]
(It is tedious but straightforward to check that each such sum is a
cocycle.)  One can see that $\massey{\alpha_{1}, \dotsc, \alpha_{n}}
\subset H^{s-(n-2)}A$, where $s$ is the sum of the degrees of the
$\alpha_{i}$'s, which means that $\massey{\alpha_{1}, \dotsc,
\alpha_{n}}$ is in the same degree as $m_{n}(\alpha_{1} \otimes \dotsb
\otimes \alpha_{n})$.  This is not a coincidence.

\begin{theorem}\label{thm-massey}
Let $A$ be a DG algebra.  Up to a sign, the higher
$A_{\infty}$-multiplications on $HA$ give Massey products.  More
precisely, if $HA$ is given an $A_{\infty}$-algebra structure by
Merkulov's construction, then for any $n \geq 3$, if $\alpha_{1},
\dotsc, \alpha_{n} \in HA$ are elements such that 
the Massey product $\massey{\alpha_{1}, \dotsc, \alpha_{n}}$ is
defined, then 
\[
(-1)^{b} m_{n}(\alpha_{1} \otimes \dotsb \otimes \alpha_{n})
\in \massey{\alpha_{1}, \dotsc, \alpha_{n}},
\]
where 
\[
b=1+\deg \alpha_{n-1} + \deg \alpha_{n-3} + \deg \alpha_{n-5} + \dotsb.
\]
\end{theorem}

The authors have been unable to find an account of this theorem in its
full generality, but for some related results, see \cite[p. 233]{Ka},
\cite[Chapter 12]{St2}, and \cite[6.3--6.4]{JL}.

Now, there are choices made in Merkulov's construction -- the choices
of the splittings \eqref{2.1.1} -- and different choices can
(depending on $A$) lead to different elements in the Massey products,
as well as different (but quasi-isomorphic) $A_{\infty}$-algebra
structures.  In any case, any choice of $A_{\infty}$-structure via
Merkulov's construction gives a ``coherent'' set of choices for an
element of each Massey product.  Of course, the
$A_{\infty}$-multiplications are also universally defined, not just
when certain products are zero.

\begin{proof}
The proof is by induction on $n$.  

The theorem discusses the situation when $n \geq 3$, but we will also
use the formula when $n=2$ in the induction: when $n=2$, we have
\[
m_{2}(\alpha_{1} \otimes \alpha_{2}) = \alpha_{1}
\alpha_{2} = \massey{\alpha_{1}, \alpha_{2}} = (-1)^{1+ \deg \alpha_{1}}
\overline{\alpha}_{1} \alpha_{2}.
\]

Now let $n=3$.  We use Merkulov's construction for the
$A_{\infty}$-algebra structure on $HA$, so we choose splittings as in
\eqref{2.1.1}, and we define the multiplication maps $m_{n}$ as in
Theorem~\ref{thm2.2}.  We use a little care when choosing the elements
$a_{ij} \in A$: we define $a_{02}$ by $G(\overline{a}_{01}a_{12}) =
a_{02}$, so $\partial(a_{02}) = \overline{a}_{01}a_{12}$ and
$G \lambda_{2} (\alpha_{1} \otimes \alpha_{2}) = (-1)^{1+\deg
\alpha_{1}} a_{02}$.  We define $a_{13}$ similarly.  Then we have
\begin{align*}
m_{3}(\alpha_{1} \otimes \alpha_{2} \otimes \alpha_{3}) &= p
\lambda_{2} (G \lambda_{1} \otimes G \lambda_{2} - G \lambda_{2}
\otimes G \lambda_{1}) (\alpha_{1} \otimes \alpha_{2} \otimes \alpha_{3}) \\
&= p\left((-1)^{1+\deg \alpha_{1}+1+\deg \alpha_{2}} a_{01} a_{13} +
(-1)^{1+1+1+\deg \alpha_{1}} a_{02} a_{23}\right) \\
&= p\left((-1)^{1+\deg \alpha_{2}} \overline{a}_{01} a_{13} +
(-1)^{1+\deg \alpha_{2}} \overline{a}_{02} a_{23} \right) \\
&= (-1)^{1+\deg \alpha_{2}} p\left( \overline{a}_{01} a_{13} +
\overline{a}_{02} a_{23} \right).
\end{align*}
(Some signs here are due to the Koszul sign convention; for example,
the map $G \lambda_{2}$ has degree 1, so $(G \lambda_{1} \otimes G
\lambda_{2}) (a_{01} \otimes a_{12} \otimes a_{23}) = (-1)^{\deg
\alpha_{1}} G\lambda_{1}(a_{01}) \otimes G\lambda_{2} (a_{12} \otimes
a_{23})$.)  

The map $p$ is the projection map from $A$ to its summand $H$.
Loosely, for any cocycle $z$, $p(z)$ is the cohomology class
represented by $z$; more precisely, $p(z)$ is the unique class in $H
\subset A$ which is cohomologous to $z$.  In the situation here, the
term in parentheses is a cocycle whose cohomology class is in
$\massey{\alpha_{1}, \alpha_{2}, \alpha_{3}}$, so we get the desired
result.

Assume that the result is true for $m_{i}$ with $i<n$.  Therefore for
all $i<j$ with $j-i<n-1$, we may choose elements $a_{i-1,j}$ by the
formula
\[
G\lambda_{j-i+1}(\alpha_{i} \otimes \dotsb \otimes \alpha_{j}) =
(-1)^{1+\deg \alpha_{j-1} + \deg \alpha_{j-3} + \dotsb} a_{i-1,j}.
\]
We write $b_{ij}$ for the exponent of $-1$ here: 
\[
b_{ij}=1+\deg \alpha_{j-1} + \deg \alpha_{j-3} + \dotsb.
\]
The last term in this sum is $\deg \alpha_{j-(2k+1)}$, where $k$ is the
maximum such that $j-(2k+1) \geq i$.  Note also for use with the Koszul
sign convention that $G\lambda_{i}$ has degree $1-i$.  Then

\newpage

\begin{align*}
m_{n}(\alpha_{1} \otimes \dotsb \otimes \alpha_{n}) &= 
p\lambda_{2} \left( \sum_{s=1}^{n-1} (-1)^{s+1} G\lambda_{s} \otimes
G\lambda_{n-s} \right) \left(a_{01} \otimes \dotsb a_{n-1,n}\right) \\
 &
\begin{aligned}
= p \Biggl( \sum_{s=1}^{n-1} & (-1)^{s+1+(1-n+s)(\deg \alpha_{1} + \dotsb + \deg
\alpha_{s})} \\
 & \ G\lambda_{s}(\alpha_{1}\otimes \dotsb \otimes \alpha_{s})
 G\lambda_{n-s}(\alpha_{s+1} \otimes \dotsb \otimes \alpha_{n}) \Biggr)
\end{aligned} \\
 & 
 =p\Biggl(\sum_{s=1}^{n-1} (-1)^{s+1+(1-n+s)(\deg \alpha_{1} + \dotsb +
 \deg \alpha_{s}) + b_{1s} + b_{s+1,n}} \ a_{0s} a_{sn}\Biggr) \\
 & 
\begin{aligned}
 =p\Biggl(\sum_{s=1}^{n-1} & (-1)^{s+1+(1-n+s)(\deg \alpha_{1} + \dotsb +
 \deg \alpha_{s}) + b_{1s} + b_{s+1,n}} \\
 & (-1)^{1-s+1+ \deg \alpha_{1} + \dotsb + \deg \alpha_{s}}
\ \overline{a}_{0s} a_{sn}\Biggr)
\end{aligned} \\
 & 
 =p\Biggl(\sum_{s=1}^{n-1} (-1)^{1+(-n+s)(\deg \alpha_{1} + \dotsb +
 \deg \alpha_{s}) + b_{1s} + b_{s+1,n}} \ \overline{a}_{0s} a_{sn}
 \Biggr) \\
 &=  p\Biggl( \sum_{s=1}^{n-1} (-1)^{b_{1n}} \ \overline{a}_{0s} a_{sn}
 \Biggr),
\end{align*}
where $b_{1n}=b$ is the sign as in the theorem: if $n-s$ is even, then the
sign is $(-1)^{1+b_{1s}+b_{s+1,n}}$, and with $n-s$ even, we have
$b_{1,n}=1+b_{1s}+b_{s+1,n}$.  If $n-s$ is odd, then the sign is 
\[
(-1)^{1+b_{s+1,n} + 1+ \deg \alpha_{s} + \deg \alpha_{s-2} + \deg
\alpha_{s-4}} = (-1)^{b_{1n}},
\]
as claimed.  As with the $n=3$ case, since the sum $\sum
\overline{a}_{0s} a_{sn}$ is a cocycle, $p$ sends it to the cohomology
class that it represents, which is an element of the Massey product
$\massey{\alpha_{1}, \dots, \alpha_{n}}$.  This finishes the proof.
\end{proof}

See Section \ref{sect6} for some examples.

\section{The bar construction and Ext}
\label{sect4}

The bar/cobar construction is one of the basic tools in homological 
algebra. Everything in this section is well-known -- see
\cite{FHT1}, for example -- but we 
need the details for the proof in the next section.

Let $A$ be a connected graded algebra and let $k$ be the trivial 
$A$-module.  Of course, the $i$-th Ext-group $\Ext^i_A(k_A,k_A)$ 
can be computed by the $i$-th cohomology of the complex 
$\Hom_A(P_A,k_A)$ where $P_A$ is any projective (or free) resolution 
of $k_A$.  Since $P_A$ is projective, $\Hom_A(P_A,k_A)$ is 
quasi-isomorphic to $\Hom_A(P_A,P_A)=\End_A(P_A)$; hence 
$\Ext^i_A(k_A,k_A)\cong H^i(\End_A(P_A))$.  Since $\End_A(P_A)$ is a 
DG algebra, the graded vector space $\Ext^*_A(k_A,k_A):=
\boplus_{i\in {\mathbb Z}} \Ext^i_A(k_A,k_A)$ has a natural algebra 
structure, and it also has an
$A_\infty$-structure by Kadeishvili's result Theorem \ref{thm2.1}. 
By \cite[Chap.2]{Ad}, the Ext-algebra of a graded algebra $A$ 
can also be computed by using the bar construction on $A$, which 
will be explained below.

First we review the shift functor. Let $(M, \partial)$ be a complex
with differential $\partial$ of degree 1, and
let $n$ be an integer. The $n$th shift of $M$, denoted by $S^n(M)$, is
defined by
\[
S^n(M)^{i}=M^{i+n}
\]
and the differential of $S^n(M)$ is 
\[
\partial_{S^n(M)}(m)=(-1)^n \partial(m)
\]
for all $m\in M$. If $f: M\to N$ is a homomorphism of degree $p$, 
then $S^n(f): S^n(M)\to S^n(N)$ is defined by the formula
\[
S^n(f)(m)=(-1)^{pn} f(m)
\]
for all $m\in S^n(M)$. The functor $S^n$ is an automorphism of 
the category of complexes. 

The following definition is essentially standard, although sign
conventions may vary; we use the conventions from 
\cite[Sect.19]{FHT1}.  Let $A$ be an augmented DG algebra with 
augmentation $\epsilon: A\to k$, viewing $k$ as a trivial DG 
algebra. Let $I$ be the kernel of $\epsilon$ and $SI$ be the 
shift of $I$. The \emph{tensor coalgebra} on $SI$ is
\[
T(SI)=k\oplus SI\oplus (SI)^{\otimes 2}\oplus (SI)^{\otimes 3} 
\oplus \cdots,
\]
where an element $Sa_1\otimes Sa_2 \otimes \cdots \otimes Sa_n$ in 
$(SI)^{\otimes n}$ is written as
\[
[a_1| a_2| \cdots |a_n]
\]
for $a_i\in I$, together with the comultiplication
\[
\Delta([a_1 |\cdots | a_n])=
\sum_{i=0}^{n} [a_1|\cdots |a_i]\otimes [a_{i+1}|\cdots |a_{n}].
\]
The degree of $[a_1|\cdots|a_n]$ is $\sum_{i=1}^n (\deg a_i -1)$.

\begin{definition}
\label{defn4.1} 
Let $(A,\partial_A)$ be an augmented DG algebra and let $I$ denote 
the augmentation ideal $\ker (A\to k)$. The {\it bar construction} 
on $A$ is the coaugmented differential graded coalgebra (DG coalgebra, 
for short) $BA$ defined as follows:

$\bullet$ As a coaugmented graded coalgebra, $BA$ is the tensor 
coalgebra $T(SI)$.

$\bullet$ The differential in $BA$ is the sum $d=d_0+d_1$ of the 
coderivations given by 
\[
d_0([a_1| \cdots | a_m])=-\sum_{i=1}^{m} (-1)^{n_i} [a_1 |\cdots| 
\partial_A(a_i)|\cdots | a_m]
\] 
and
\begin{gather*}
d_1([a_1])=0
\\
d_1([a_1 |\cdots |a_m])=\sum_{i=2}^{m} (-1)^{n_i} [a_1|\cdots 
|a_{i-1}a_i| \cdots |a_m]
\end{gather*}
where $n_i=\sum_{j<i} (-1+\deg a_j)=\sum_{j<i} \deg [a_j]$. 
\end{definition}

The cobar construction $\Omega C$ on a coaugmented DG coalgebra 
$C$ is defined dually 
\cite[Sect.19]{FHT1}. We omit the definition since it is used only 
in two places, one of which is between Lemma \ref{lem5.3}
and Lemma \ref{lem5.4}, and the other is in Lemma \ref{lem5.5}.

In the rest of this section we assume that 
$A$ is an augmented associative algebra.
In this case $SI$ is concentrated in degree $-1$; hence the degree 
of $[a_1|\cdots |a_{m}]$ is $-m$.  This means that
the bar construction $BA$ is graded by the {\it negative} of tensor 
length. The degree of the differential $d$ is $1$. We may think 
of $BA$ as a complex with $(-i)$th term equal to $I^{\otimes i}$,
the differential $d$ mapping
$I^{\otimes i}$ to $I^{\otimes i-1}$. If $A$ has an Adams grading, 
denoted $\adeg$, then $BA$ has a bigrading that is defined by
\[
\deg\; [a_1|\cdots |a_{m}]=(-m, \sum_i \adeg a_i).
\]
The second component is the Adams degree of $[a_1|\cdots |a_{m}]$.  

The bar construction on the left $A$-module $A$, denoted by $B(A,A)$, 
is constructed as follows. As a complex $B(A,A)=BA\otimes A$ with 
$(-i)$th term equal to $I^{\otimes i}\otimes A$. We use 
\[
[a_1|\cdots|a_m]x
\]
to denote an element in $I^{\otimes i}\otimes A$ where $x\in A$ and 
$a_i\in I$. The degree of $[a_1|\cdots|a_m] x$ is $-m$. The 
differential on $B(A,A)$ is defined by
\[
d(x)=0 \ \ \text{($m=0$ case)},
\]
and
\[
d([a_1 |\cdots | a_m]x)=\sum_{i=2}^{m} (-1)^{i-1} [a_1|\cdots |a_{i-1}a_i|
\cdots |a_m]x+(-1)^m [a_1|\cdots |a_{m-1}]a_m x. 
\]
Then $B(A,A)$ is a complex of free right $A$-modules.
One basic property is that the augmentations of $BA$ and $A$ make it
into a free resolution of $k_A$,
\begin{equation}
\label{4.1.1}B(A,A)\to k_A\to 0
\end{equation}
(see \cite[19.2]{FHT1} and \cite[Chap.2]{Ad}). 

\begin{remark}
\label{rem4.2}
In the next section we use the tensor $\otimes$ notation instead 
of the bar $|$ notation, which seems more natural when we concentrate 
on each term of the bar construction. 
\end{remark}

We now assume that with respect to the Adams grading, $A$ is connected
graded and finite-dimensional in each degree.
Then $B(A,A)$ is bigraded with Adams grading on the second 
component, and the differential of $B(A,A)$ preserves the Adams grading.
Let $B^{\#}A$ be the graded $k$-linear dual of the coalgebra $B A$. 
Since $BA$ is locally finite, $B^{\#}A$ is a 
locally finite bigraded algebra.  With respect to the Adams grading, 
$B^{\#}A$ is {\it negatively} connected graded. The DG algebra 
$\End_A(B(A,A)_A)$ is bigraded too, but not Adams connected. Since 
$B(A,A)$ is a left differential graded comodule over $BA$, it has a 
left differential graded module structure over $B^{\#}A$, which is 
compatible with the right $A$-module structure. By an 
idea  similar to \cite[Ex. 4, p.\ 272]{FHT1} (also see \cite{LP2}) 
one can show  that the natural map $B^{\#}A\to \End_A(B(A,A)_A)$ is 
a quasi-isomorphism of DG algebras. 

Define the \emph{Koszul dual} of a connected graded ring $A$ to be the 
DG algebra $\End_A(P_A)$, where $P_A$ is any free resolution of 
$k_A$.  By the following lemma,
this definition makes sense up to quasi-isomorphism in the category of 
$A_\infty$-algebras.

\begin{lemma}
\label{lem4.3} 
Let $A$ be a connected graded algebra which is finite-dimensional in
each degree, and let $P_A$ and $Q_A$
be two free resolutions of $k_A$.
\begin{enumerate}
\item $\End_A(P_A)$ is quasi-isomorphic to $\End_A(Q_A)$ as 
$A_\infty$-algebras.
\item $\End_A(P_A)$ is quasi-isomorphic to $B^{\#}A$ as $A_\infty$-algebras.
\end{enumerate}
\end{lemma}

\begin{proof}
(a) We may assume that $Q_A$ is a minimal free resolution of $k_A$. Then 
$P_A=Q_A\oplus I_A$ where $I_A$ is another complex of free modules
such that $HI_A=0$ \cite[10.1.3 and 10.3.4]{AFH}. 
In this case $D:=\End_A(Q_A)$ is a sub-DG algebra of 
$E:=\End_A(P_A)$ such that $D=eEe$ where $e$ is the projection onto
$Q_A$. Let $C=(1-e)E+E(1-e)$. Then 
\[ 
C=\Hom_A(I_A,Q_A)+\Hom_A(Q_A,I_A)+\Hom_A(I_A,I_A),
\]
and $HC=0$. By Lemma \ref{lem2.5}, $D$ and $E$ are quasi-isomorphic.

(b) Since $B(A,A)$ is a free resolution of $k_A$, then part (a) says
that $\End_A(P_A)$ is 
quasi-isomorphic to $\End_A(B(A,A)_A)$. The assertion follows
from the fact that $\End_A(B(A,A))$ is quasi-isomorphic to $B^{\#}A$
\cite[Ex.\ 4, p.\ 272]{FHT1}.
\end{proof}

So we may think of the bigraded DG algebra $B^{\#}A$ as the Koszul dual of
$A$.  This viewpoint of Koszul duality is also taken by Keller in
\cite{Ke1}.  By results in \cite{LP2}, we can define the Koszul dual
of any connected graded (or augmented) $A_\infty$-algebra, and the
double Koszul dual is quasi-isomorphic to the original
$A_\infty$-algebra.

The classical Ext-algebra $\Ext^*_A(k_A,k_A)$ is the cohomology ring
of $\End_A(P_A)$, where $P_A$ is any free resolution of $k_A$. The
above lemma demonstrates the familiar fact that this is independent of
the choice of $P_A$.  Since $E:=\End_A(P_A)$ is a DG algebra, by Proposition
\ref{prop2.3}, $\Ext^*_A(k_A,k_A)=HE$ has a natural
$A_\infty$-structure, which is called an {\it $A_\infty$-Ext-algebra}
of $A$. By abuse of notation we use $\Ext^*_A(k_A,k_A)$ to denote an
$A_\infty$-Ext-algebra.

\section{$A_\infty$-structure on Ext-algebras} 
\label{sect5}

In this section we consider the multiplications on an
$A_\infty$-Ext-algebra of a connected graded algebra, and finally give
proofs of Theorem \ref{thmA} and Proposition~\ref{prop1.2}.  Consider
a connected graded algebra
\[
A=k\oplus A_1\oplus A_2\oplus \cdots,
\]
which is viewed as an $A_\infty$-algebra concentrated in degree 0,
with the grading on $A$ being the Adams grading. Let $Q\subset A$ be a
minimal graded vector space which generates $A$. Then $Q\cong
\fm/\fm^2$ where $\fm:=A_{\geq 1}$ is the unique maximal graded ideal
of $A$. Following Milnor and Moore \cite[3.7]{MM}, we call the
elements of $\fm / \fm^{2}$ the \emph{indecomposables} of $A$, and by
abuse of notation, we also call the elements of $Q$ indecomposables.
Let $R\subset \free{Q}$ be a minimal graded vector space which
generates the relations of $A$ ($R$ is not unique). Then $A\cong
\free{Q}/(R)$ where $(R)$ is the ideal generated by $R$, and the start
of a minimal graded free resolution of the trivial right $A$-module
$k_A$ is
\begin{equation}
\label{5.0.1}
\cdots \to R\otimes A\to Q\otimes A \to A\to k\to 0.
\end{equation}

\begin{lemma}
\label{lem5.1} 
Let $A$ be a connected graded algebra. 
Then there are natural isomorphisms of graded vector spaces 
\[
\Ext^1_{A}(k_A, k_A)\cong Q^{\#} = \boplus Q_{s}^{\#}
\quad \text{and}\quad \Ext^2_{A}(k_A,k_A)\cong R^{\#} = \boplus
R_{s}^{\#}.
\]
\end{lemma}

\begin{proof} 
This follows from the minimal free resolution \eqref{5.0.1}.
\end{proof}

In the rest of the section, we assume that $A$ is Adams locally
finite: each $A_{i}$ is finite-dimensional.  Let $E$ be the
$A_\infty$-Ext-algebra $\Ext^*_{A}(k_A,k_A)$.  We would like to
describe the $A_\infty$-structure on $E$ by using Merkulov's
construction.

We first fix some notation. For each Adams degree $s$, we choose a
vector space splitting $A_{s} = Q_{s} \oplus D_{s}$; the elements in
$Q_{s}$ are indecomposable, while those in $D_{s}$ are
``decomposable'' in terms of the indecomposables of degree less than
$s$.
More precisely, we define $Q_{s}$ and $D_{s}$ inductively: we start by
setting $Q_{1}=A_{1}$ and $D_{1}=0$.  Now assume that $Q_{i}$ and
$D_{i}$ have been defined for $i < s$; then we set $D_{s}$ to be the
image in $A_{s}$ of the multiplication map
\[
\mu_{s} : \boplus_{1 \leq i \leq s-1} Q_{i} \otimes A_{s-i}
\rightarrow A_{s}.
\]
We choose $Q_{s}$ to be a vector space complement of $D_{s}$.  

Now for each $s \geq 2$, the multiplication map
\[
\mu_{s} : \boplus_{1 \leq i \leq s} Q_{i} \otimes A_{s-i} \rightarrow A_{s}
\]
is onto.  Define the $k$-linear map $\xi_{s}: A_{s} \rightarrow
\boplus_{1 \leq i \leq s} Q_{i} \otimes A_{s-i}$ so that the composition 
\begin{equation}
\label{5.1.1}
A_s \xrightarrow{\xi_s} \boplus_{1 \leq i \leq s} Q_{i} \otimes A_{s-i}
\xrightarrow{-\mu_s} A_s
\end{equation}
is the identity map of $A_s$.  Further, we choose $\xi_{s}$ so that
with respect to the direct sum decomposition $A_{s} = Q_{s} \oplus
D_{s}$, we have 
\begin{equation}\label{5.1.2}
\im (\xi_{s}\vert_{Q_{s}}) = Q_{s} \otimes A_{0}, \quad 
\im (\xi_{s}\vert_{D_{s}}) \subset \boplus_{1 \leq i \leq s-1}
Q_{i} \otimes A_{s-i}.
\end{equation}
(The second of these holds for any choice of $\xi_{s}$; the first need
not.)  Define $\xi_1=\theta_1=id_{A_1}$, and inductively set
$\theta_s=\sum_{i+j=s} (id_{Q_{i}}\otimes \theta_{j}) \circ \xi_s$;
that is, $\theta_{s}$ is the composition
\begin{gather*}
A_s \xrightarrow{\xi_s} \boplus_{i+j =s} Q_{i} \otimes A_{j}
\xrightarrow{\sum id_{Q_{i}}\otimes \xi_{j}}
\boplus_{i+k+l=s} Q_{i} \otimes Q_{k} \otimes A_{l} 
\xrightarrow{\sum id_{Q_{i}} \otimes id_{Q_{k}} \otimes \xi_{l}}
\cdots \\
\longrightarrow
\boplus_{n \geq 1} \boplus_{i_{1}+\dotsb +i_{n}=s} Q_{i_{1}} \otimes
\dotsb \otimes Q_{i_{n}}.
\end{gather*}
Here, the subscripts on the $Q$'s are positive, while those
on the $A$'s are non-negative.

Let $R=\boplus_{s\geq 2} R_{s} \subset \free{Q}$ be a minimal graded
vector space of the relations of $A$.  Note that with respect to
tensor length, the elements of $R$ need not be homogeneous, but they are
homogeneous with respect to the Adams grading -- the grading induced
by that on $Q$.

Let $\free{Q}_{s}$ denote the part of $\free{Q}$ in Adams degree $s$;
thus
\[
\free{Q}_{s} = \boplus_{n \geq 1} \boplus_{i_{1}+\dotsb
+i_{n}=s} Q_{i_{1}} \otimes \dotsb \otimes Q_{i_{n}}.
\]
We write $\mu$ for the map $\mu : \free{Q} \rightarrow A \cong
\free{Q}/(R)$.

\begin{lemma}\label{lem5.2}
For each $s$, $R_{s}$ may be chosen so that
\[
R_s\subset \boplus_{1 \leq i \leq s-1}
Q_{i}\otimes \theta_{s-i}(A_{s-i}) \subset \free{Q}_{s}.
\]
Hence $R_{s}$ may also be viewed as a subspace of $\boplus_{1 \leq i
\leq s-1} Q_{i} \otimes A_{s-i}$, via the composite
\[ 
R_s \hookrightarrow 
\boplus_{1 \leq i \leq s-1} Q_{i}\otimes \theta_{s-i}(A_{s-i})
\xrightarrow{\sum id_{Q_{i}} \otimes \mu} 
\boplus_{1 \leq i \leq s-1} Q_{i} \otimes A_{s-i}.
\]
\end{lemma}

\begin{proof}
Let $(R)_{s-i}$ be the degree $s-i$ part of the ideal $(R)$. Then it
is generated by the relations of degree at most $s-i$, and we have a
decomposition
\[
\free{Q}_{s-i} = \theta_{s-i}(A_{s-i}) \oplus \ker \mu=
\theta_{s-i}(A_{s-i})\oplus (R)_{s-i},
\]
where $\mu: \free{Q}_{s-i}\to A_{s-i}$ is multiplication.  
Hence we have 
\[ 
\free{Q}_{s} = Q_{s} \oplus \boplus_{1 \leq i \leq s-1} 
\bigl[ \bigl( Q_{i} \otimes \theta_{s-i} (A_{s-i})\bigr) 
\oplus \bigl(Q_{i} \otimes (R)_{s-i}\bigr) \bigr].
\]
Any relation $r\in R_s$ has no summands in $Q_{s}$, and hence is a sum
of $r'\in \boplus Q_{i} \otimes \theta_{s-i} (A_{s-i})$ and $r''\in
\boplus Q_{i} \otimes (R)_{s-i}$. Modulo the relations of degree less
than $n$, we may assume $r''=0$. Hence the first part of the lemma is
proved.

The map $\theta_{s-i}: A_{s-i} \rightarrow \free{Q}_{s-i}$ is an
inclusion, and up to a sign its left inverse is the multiplication map
$\mu: \free{Q}_{s-i} \rightarrow A_{s-i}$.  Once $R_{s}$ has been
chosen to be a subspace of $\boplus Q_{i} \otimes \theta_{s-i}
(A_{s-i})$, composing with $\mu$ takes it injectively to $\boplus
Q_{i} \otimes A_{s-i}$.
\end{proof}

The minimal resolution \eqref{5.0.1} is a direct summand of any other
resolution, and in particular it is a summand of the bar
resolution \eqref{4.1.1}
\[ 
\cdots \to \fm^{\otimes 2}\otimes A\to \fm\otimes A \to A\to k\to 0. 
\]
We have made several choices up to this point: choosing the splittings
$A_{s}=Q_{s} \oplus D_{s}$, and now choosing $R_s$ as in the lemma,
so that $R_s\subset \boplus_{1 \leq i \leq s-1} Q_{i} \otimes A_{s-i}
\subset \fm\otimes \fm$.  These choices give a choice for this
splitting of resolutions, at least in low degrees.

Since $A$ is concentrated in degree 0, the grading on the differential 
graded coalgebra $T(S\fm)$ is by the negative of the wordlength, namely, 
$(T(S\fm))^{-i}=\fm^{\otimes i}$. The differential $d=(d^i)$ of the bar 
construction $T(S\fm)$ is induced by the multiplication 
$\fm \otimes \fm \to \fm$ in $A$. For example, 
\[
d^{-1}([a_1])=0 \quad \text{and}\quad d^{-2}([a_1|a_2])=(-1)^{-1}[a_1 a_2]
\]
for all $a_1, a_2\in \fm$. 
There is a natural decomposition of $\fm$ with respect to the Adams
grading,
\[
\fm=A_1\oplus A_2\oplus A_3\oplus \cdots,
\]
which gives rise to a decomposition of $\fm\otimes \fm$ with respect to the 
Adams grading:
\[
\fm\otimes \fm=(A_1\otimes A_1)\oplus (A_1\otimes A_2\oplus
A_2\otimes A_1) \oplus \cdots.
\]
As mentioned above, we are viewing $R_s$ as being a subspace of
$\boplus Q_i\otimes A_{s-i}$. 

\begin{lemma}
\label{lem5.3} 
Let $W_{s}$ be the Adams degree $s$ part of $\fm \otimes \fm$;
that is, let 
\[
W_{s}=\boplus_{1 \leq i \leq s-1} A_i\otimes A_{s-i}.
\]
Then there are decompositions of vector spaces
\begin{gather*}
W_{s}=\im (d^{-3}_{s})\oplus R_s\oplus \xi_s(D_s), \\
\ker (d^{-2}_{s})=\im (d^{-3}_{s}) \oplus R_s,
\end{gather*}
where $R_s$ and $\xi_s(D_s)$ are subspaces of 
\[
\boplus_{1 \leq i \leq s-1} Q_{i} \otimes A_{s-i} \subset 
\boplus_{1 \leq i \leq s-1} A_{i} \otimes A_{s-i} = W_{s}.
\]
\end{lemma}

\begin{proof} It is clear that the injection 
\[
D_s \xrightarrow{\xi_s} \boplus_{1 \leq i \leq s-1} Q_{i} \otimes
A_{s-i} \longrightarrow W_{s}
\]
defines a projection from $W_s$ to $D_s$. Since $d^{-2}_{s}:
W_s\to D_s$ is a surjection, we have a decomposition
\[
W_s=\ker (d^{-2}_{s})\oplus \xi_s(D_s).
\]
Since $R^{\#}\cong \Ext^2_{A}(k_A,k_A)=H^2((T(S\fm))^{\#})$ by Lemma
\ref{lem5.1}, there is a decomposition $\ker (d^{-2}_{s})=\im
(d^{-3}_{s}) \oplus R_s$.  Hence the assertion follows.
\end{proof}

Since $A$ is Adams locally finite, $(\fm^{\otimes n})^{\#}\cong
(\fm^{\#})^{\otimes n}$ for all $n$. Let $\Omega A^{\#}$ be the cobar 
construction on the DG coalgebra $A^{\#}$. Via the isomorphisms
\[
B^{\#}A=(T(S\fm))^{\#}\cong T((S\fm)^{\#})\cong T(S^{-1}\fm^{\#})
=\Omega A^{\#},
\]
we identify $B^{\#}A=(T(S\fm))^{\#}$ with $\Omega A^{\#}=T(S^{-1}\fm^{\#})$.
The differential $\partial$ on $B^{\#}A$ is defined by
\[
\partial(f)=-(-1)^{\deg f} f\circ d
\]
for all $f\in T(S^{-1}\fm^{\#})$ .

We now study the first two nonzero differential maps of $\Omega A^{\#}$,
\[
\partial^1: \fm^{\#}\to (\fm^{\#})^{\otimes 2}\qquad\text{and}\qquad 
\partial^2: (\fm^{\#})^{\otimes 2}\to  (\fm^{\#})^{\otimes 3}.
\] 
For all $s$ and $n$, let
\[
T^n=(\fm^{\#})^{\otimes n}
\]
and
\[
T^n_{-s}=\boplus_{i_1+\cdots +i_n=s} A^{\#}_{i_1}\otimes \cdots \otimes
A^{\#}_{i_n}.
\]
Since we are working with $\fm^{\#}$, all subscripts here and in what
follows are positive.  Fix Adams degree $-s$, and consider
\begin{gather*}
\partial^{1}_{-s}:A_s^{\#}\to \boplus_{i+j=s} A_i^{\#}\otimes A_j^{\#}, \\
\partial^{2}_{-s}:\boplus_{i+j=s} A_i^{\#}\otimes A_j^{\#}\to 
\boplus_{i_1+i_2+i_3=s}A_{i_1}^{\#}\otimes A_{i_2}^{\#}\otimes A_{i_3}^{\#}.
\end{gather*}
The decomposition \eqref{2.1.1} for $T^1_{-s}$ is
\[
B^{1}_{-s} = 0, \quad H^{1}_{-s} = Q^{\#}_{s}, \quad L^{1}_{-s} =
D^{\#}_{s}.
\]
for all $s\geq 1$. The decomposition \eqref{2.1.1} for $T^2_{-s}$ is
given in the following lemma.

\begin{lemma} 
\label{lem5.4} Fix $s \geq 2$.  With notation as above, we have the
following.
\begin{enumerate}
\item Define the duals of subspaces by using the decompositions given in 
Lemma \ref{lem5.3}. Then $\im  \partial^1_{-s}=(\xi_s(A_s))^{\#}$ and 
$\ker \partial^2_{-s} =(\xi_s(A_s))^{\#}\oplus R_s^{\#}.$
\item 
The decomposition \eqref{2.1.1} for $T^{2}_{-s}$ can be chosen to be
\[
T^2_{-s}=B^2_{-s}\oplus H^2_{-s} \oplus L^2_{-s}=
(\xi_s(A_s))^{\#}\oplus R_s^{\#}\oplus (\im d^{-3}_{s})^{\#}.
\] 
The projections onto $R_s^{\#}$ and $(\xi_s(A_s))^{\#}$ kill
$\boplus_{2 \leq i \leq s-1} D_i^{\#}\otimes A_{s-i}^{\#}$.
\item Let $G$ be the homotopy defined in Merkulov's construction for the
DG algebra $T(\fm^{\#})$. Then we may choose $G^{2}_{-s}$ to be equal to
$-(\xi_s)^{\#}$, restricted to $T^2_{-s}$.
\end{enumerate}
\end{lemma} 

\begin{proof} (a) This follows from Lemma \ref{lem5.3} and a linear 
algebra argument. 

(b) This follows from Lemma \ref{lem5.3}, part (a), and the fact that
$R_s$ and $\xi_s(A_s)$ are subspaces of $\boplus_{1 \leq i \leq s-1}
Q_{i} \otimes A_{s-i}$.

(c) Let $\xi_s$ also denote the map $D_s\to \boplus_{1 \leq i \leq
s-1} Q_{i} \otimes A_{s-i} \to W_s$.  Since $d^{-2}_s=-\mu_s$, the
composite $d^{-2}_s \circ \xi_s: D_s\to W_s\to D_s$ is the identity
map (see \eqref{5.1.1}).  Since $ (d^{-2}_s \circ
\xi_s)^{\#}=(\xi_s)^{\#} \circ (d^{-2}_s)^{\#}$,
\[
(\xi_s)^{\#}\circ (d^{-2}_s)^{\#}: \ D_s^{\#}\to T^2_{-s}\to D_s^{\#}
\]
is the identity map of $D_s^{\#}$, which is the summand $L^{1}_{-s}$
of $T^1$.  Since $(d^{-2}_s)^{\#}=-\partial^1_{-s}$, we may choose the
homotopy $G$ to be $-(\xi_s)^{\#}$ when restricted to $T^2_{-s}$.
\end{proof}

Now we start to construct the higher $A_{\infty}$-multiplication maps
on $\Ext^{*}_{A}(k_{A},k_{A})$, using Merkulov's construction from
Section \ref{sect2}.  Lemma \ref{lem5.4} tells us what the homotopy
$G$ is.  The maps $\lambda_{n} : (T(S^{-1}\fm^{\#}))^{\otimes n}
\rightarrow T(S^{-1}\fm^{\#})$ are defined as in Section 2; in
particular, recall that we formally set $G\lambda_1=-id_T$, and
$\lambda_2$ is the multiplication of $T(S^{-1}\fm^{\#})$.

Recall that $T(S^{-1}\fm^{\#})$ is a free (or tensor) DG algebra
generated by $S^{-1}\fm^{\#}$.  To distinguish among the various
tensor products occurring here, we use $\botimes$ when tensoring
factors of $T(S^{-1}\fm^{\#})$ together; in particular, we write
$\lambda_{2}$ as
\[
\lambda_{2} : (S^{-1}\fm^{\#})^{\otimes n} \botimes
(S^{-1}\fm^{\#})^{\otimes m} \rightarrow (S^{-1}\fm^{\#})^{\otimes (n+m)}.
\]
Then for $a_1\otimes\cdots\otimes a_n\in 
(S^{-1}\fm^{\#})^{\otimes n}$ and $b_1\otimes \cdots \otimes b_m
\in (S^{-1}\fm^{\#})^{\otimes m}$ we have
\begin{equation}
\label{A.12.1}
\lambda_2((a_1\otimes\cdots\otimes a_n)\botimes (b_1\otimes \cdots
\otimes b_m))=
a_1\otimes \cdots\otimes a_n\otimes b_1\otimes \cdots\otimes b_m.
\end{equation}
By the above formula, we see that $\lambda_2$ changes $\botimes$ to $\otimes$,
so it is like the identity map. 

\begin{lemma} 
\label{lem5.5} 
Let $E^1=\Ext^1_{A}(k_A,k_A)=Q^{\#}$. Fix $n \geq 2$ and $s\geq 2$.
\begin{enumerate}
\item When restricted to $(E^1)^{\otimes n}$ in Adams degree $s$, the
map $\lambda_n$ has image in 
\[
T^2_{-s} =\boplus_{q+r=s} A_q^{\#}\otimes A_r^{\#}.
\]
Hence the image of $G \lambda_{n}$ is in $D_{s}^{\#}$.
\item When restricted to $(E^1)^{\otimes n}$ in Adams degree $s$, the
map $-G\lambda_n$ is the $k$-linear dual of the composite 
\[
A_s  \xrightarrow{\theta_{s}} 
\boplus_{m \geq 1} \boplus_{i_{1}+\dotsb +i_{m}=s} Q_{i_{1}} \otimes \dotsb
\otimes Q_{i_{m}} 
 \rightarrow \boplus_{i_{1}+\dotsb +i_{n}=s} Q_{i_{1}} \otimes \dotsb
\otimes Q_{i_{n}}.
\]
\item When restricted to $(E^1)^{\otimes n}$ in Adams degree $s$, the
map $m_n=Pr_H \lambda_n$ is the $k$-linear dual of the canonical map
\[
R_s\longrightarrow 
\boplus_{1 \leq i \leq s-1} Q_{i} \otimes A_{s-i}
\xrightarrow{\sum id \otimes \theta_{s-i}} 
\boplus_{i_{1}+\dotsb +i_{n}=s} Q_{i_{1}} \otimes \dotsb \otimes Q_{i_{n}}.
\]
\end{enumerate}
\end{lemma} 

\begin{proof} We use induction on $n$. 

(a) By definition,
\[
\lambda_n=\lambda_2 \sum_{\substack{i+j=n \\ i,j>0}} (-1)^{i+1}
G\lambda_i \botimes G \lambda_j.
\]
For $n=2$, the claim follows from \eqref{A.12.1}. Now assume $n>2$ and
consider $\lambda_{n}$ applied to $b_{1} \otimes \dotsb \otimes
b_{n}$, with $b_{m} \in Q_{i_{m}}^{\#}$ for each $m$.
By Lemma~\ref{lem5.4}(c), when restricted to $T^2_{-r}$ for any $r$,
$-G$ is dual to the map $\xi_{r}: A_{r} \rightarrow \boplus Q_{i}
\otimes A_{r-i} \subset \boplus A_{i} \otimes A_{r-i}$, so by
induction, for each $m < n$,
\[
G \lambda_{m} (b_{1} \otimes \dotsb \otimes b_{m}) \in G(
\boplus_{i+j=i_{1}+\dotsb+i_{m}} A_{i}^{\#} \otimes A_{j}^{\#}) 
\subset A^{\#}_{i_{1}+\dotsb +i_{m}},
\]
and similarly for $G \lambda_{n-m}(b_{m+1} \otimes \dotsb \otimes
b_{n})$.  Thus the first statement follows.  The second statement
follows from the assumption \eqref{5.1.2} about the map $\xi_{n}$.

(b) When $n=2$, $\theta_2=\xi_2$, and the claim follows from Lemma
\ref{lem5.4}(c).  Now we assume $n>2$. When restricted to
$(E^1)^{\otimes n}$ in Adams degree $s$, part (a) says that if $i>1$,
then the image of $G\lambda_i\botimes G\lambda_j$ is in
$\boplus_{q+r=s} D_{q}^{\#} \botimes D_{r}^{\#} \subset
\boplus_{q+r=s} D_q^{\#}\botimes A_r^{\#}$. By Lemma
\ref{lem5.4}(b,c),
\[
G\lambda_2(D_q^{\#}\botimes A_r^{\#})=G(D_q^{\#}\otimes A_r^{\#})=
-(\xi_n)^{\#}(D_q^{\#}\otimes A_r^{\#})=0.
\]
Therefore, when restricted to $((E^1)^{\otimes n})_{s}$, we have 
\[
G\lambda_n=G\lambda_2[(-1)^{2} (-id)\botimes G\lambda_{n-1}]=
-G\lambda_2(id\botimes G\lambda_{n-1}).
\]
By induction on $n$, in any Adams degree $r$, $G\lambda_{n-1}$ is
$-(\theta_{r})^{\#}$ composed with projection to a summand, and by Lemma 
\ref{lem5.4}(c) we see that in Adams degree $q$, we have
$G=-(\xi_q)^{\#}$.  Hence when applied to $b_{1} \otimes b_{2} \otimes
\dotsb \otimes b_{n}$ with $b_{1}$ in Adams degree $i$ and $b_{2}
\otimes \dotsb \otimes b_{n}$ in Adams degree $s-i$, we have
\[
G\lambda_n = (\xi_s)^{\#} \lambda_2 (id_{E^1_{i}} \botimes
-(\theta_{s-i})^{\#}) = -((id_{Q_{i}}\otimes \theta_{s-i})\circ
\xi_s)^{\#},
\]
and thus $G \lambda_n$ is exactly $-\theta_{s}^{\#}$ followed by
projection onto the tensor length $n$ summand 
\[
\boplus_{i_{1}+\dotsb +i_{n}=s} Q_{i_{1}} \otimes \dotsb \otimes
Q_{i_{n}}.
\]

(c) Since we assume that $R_s$ is a subspace of $\boplus_{1 \leq
i \leq s-1} Q_{i} \otimes A_{s-i}$, the dual of the inclusion
\[
R_s\to \boplus_{1 \leq i \leq s-1} Q_{i} \otimes A_{s-i}
\]
is $Pr_H$ restricted to $\left( \boplus_{1 \leq i \leq s-1} Q_{i}
\otimes A_{s-i} \right)^{\#}$. Hence the dual of
\[
R_s\to \boplus_{1 \leq i \leq s-1} Q_{i} \otimes A_{s-i} \rightarrow
(Q \otimes Q^{\otimes (n-1)})_{s}
\]
is equal to $Pr_H\circ (\sum id\otimes \theta_{s-i})^{\#}$.

By Lemma \ref{lem5.4}(b), $Pr_H$ is zero when applied to 
$D_q^{\#}\otimes A_r^{\#}$ for all $q$. By (b),
\[
Pr_H \lambda_n =Pr_H \lambda_2 \left( -\sum_{q+r=s}
id_{Q^{\#}_q}\botimes G\lambda_{r} \right) =Pr_H \left( \sum_{q+r=s}
id_{E^1_q}\otimes (\theta_r)^{\#}\right),
\]
which is the desired map.
\end{proof} 

\begin{proof}[Proof of Theorem \ref{thmA}] 
First of all by Lemma \ref{lem5.2}, we may assume that 
\[
R_{s}\subset \boplus_{1 \leq i \leq s-1} Q_{i} \otimes A_{s-i}.
\]
Then we appeal to Lemmas \ref{lem5.3}, \ref{lem5.4} and
\ref{lem5.5}. The canonical map in Lemma \ref{lem5.5}(c) is just the
inclusion, and so the assertion holds.
\end{proof}

\begin{proof}[Proof of Corollary \ref{corB}]
Note that under the assumptions in the corollary, the vector space of
indecomposables is canonically isomorphic to $A_{1}$.  Thus various
parts of Theorem \ref{thmA} simplify; for example, $E^{1}$ is
isomorphic to $A_{1}^{\#}$, and hence is concentrated in Adams degree
$-1$.
\end{proof}

The following corollary is immediate.

\begin{corollary} 
\label{cor5.6} Let $A$ and $E$ be as in Theorem \ref{thmA}. 
\begin{enumerate}
\item The algebra $A$ is determined by the 
maps $m_n$ restricted to $(E^1)^{\otimes n}$ for all $n$.
\item The $A_\infty$-structure of $E$ is determined up to 
quasi-isomorphism by the 
maps $m_n$ restricted to $(E^1)^{\otimes n}$ for all $n$.
\end{enumerate}
\end{corollary}

\begin{proof}
(a) By Theorem \ref{thmA}, the map $R\to \free{Q}$ can be recovered
from $m_n$ restricted to $(E^1)^{\otimes n}$. Hence the structure of
$A$ is determined.

(b) After $A$ is recovered, the $A_\infty$-structure of $E$ is determined 
by $A$. Therefore the structure of $E$ is determined by the 
restriction of $m_n$ on $(E^1)^{\otimes n}$, up to quasi-isomorphism.
\end{proof}

\begin{proof}[Proof of Proposition \ref{prop1.2}]
By the construction given above, it is clear that if the grading group
for the Adams grading is $\Z\oplus G$ for some abelian group $G$, then
all of the maps including $m_n$ preserve the $G$-grading.  The
assertion follows.
\end{proof}

\section{Examples}
\label{sect6}

Several examples of $A_\infty$-algebras $E$ are given in
\cite[Exs. 3.5, 3.7, 13.4 and 13.5]{LP1}. We conclude this paper with
a few more examples.

\begin{example}\label{ex5.7}
Fix a field $k$ and consider the free algebra $B=k\langle x_{1}, x_{2}
\rangle$ on two generators, each in Adams degree 1.  Fix an integer $q
\geq 2$, and let $f(x_{1},x_{2})$ be an element in Adams degree $q$,
and let $A=B/(f)$.  Then the minimal resolution \eqref{5.0.1} for $k_{A}$
has the form
\[
\dotsb \rightarrow Ar \xrightarrow{i} Ae_{1} \oplus Ae_{2} \rightarrow
A \rightarrow k_{A} \rightarrow 0,
\]
where the generator $e_{i}$ corresponds to $x_{i}$, and $r$
corresponds to $f(x_{1},x_{2})$.  Order the monomials in $B$
left-lexicographically, setting $x_{1}<x_{2}$, and assume that with
respect to this ordering, $f(x_{1},x_{2})$ has leading term
$x_{2}^{i}x_{1}^{q-i}$ with $0<i<q$.  Then one can show that the map
$i: Ar \rightarrow Ae_{1} \oplus Ae_{2}$ is injective, so
\[
\Ext^{s}_{A}(k_{A}, k_{A}) = \begin{cases}
k & s=0, \\
k(-1) \oplus k(-1) & s=1, \\
k(-q) & s=2, \\
0 & \text{else}.
\end{cases}
\]
$\Ext^{1}$ is dual to $A_{1}$, and we choose $(y_{1}, y_{2})$ to be
the dual basis to $(x_{1}, x_{2})$.  We write $z$ for the generator of
$\Ext^{2}$ dual to $f$. For degree reasons, the $A_{\infty}$-algebra
structure on Ext has the property that $m_{n}=0$ unless $n=q$.  By
Theorem \ref{thmA}, the map $m_{q}$ is ``dual to the relations'':
\[
m_{q}(y_{i_{1}} \otimes \dotsb \otimes y_{i_{q}}) = \alpha z \quad \text{if
$\alpha x_{i_{1}} \dotsb x_{i_{q}}$ is a summand in $f(x_{1}, x_{2})$}.
\]
So for example, if $q>2$, then as an associative algebra,
$\Ext_{A}^{*}(k_{A}, k_{A})$ has trivial multiplication no matter what
$f$ is, so one cannot recover $A$ from the ordinary algebra
structure.  One can recover $A$ from the $A_{\infty}$-algebra
structure, though.
\end{example}

\begin{example}
Let $A=k[x_{2}, x_{3}]/(x_{3}^{2}-x_{2}^{3})$, graded by giving each
$x_{i}$ Adams degree $i$.  The graded vector space $Q$ has two nonzero
graded pieces: $Q_{i}$ is spanned by $x_{i}$ when $i=2,3$.  Let
$(b_{2}, b_{3})$ be the graded basis for $Q^{\#}$ which is dual to
$(x_{2}, x_{3})$.  The space $R$ of relations has two graded pieces
also: there is the degree 5 relation $r_{5}=x_{2}x_{3}-x_{3}x_{2}$,
and the degree 6 relation $r_{6}=x_{3}^{2}-x_{2}^{3}$.  Let
$(s_{5},s_{6})$ be the graded basis for $R^{\#}$ which is dual to the
basis $(r_{5},r_{6})$.  Thus in low dimensions, the Ext algebra is
given by
\[
\Ext^{n}_{A}(k_{A}, k_{A}) = \begin{cases}
k & n=0, \\
k(-2) \oplus k(-3) & n=1, \\
k(-5) \oplus k(-6) & n=2. \\
\end{cases}
\]
Indeed, by viewing $A$ as a subalgebra of $k[y]$ (with $A
\hookrightarrow k[y]$ defined by $x_{i} \mapsto y^{i}$), one can
construct a minimal resolution for $k_{A}$ to find that if $n>0$,
\[
\Ext^{n}_{A}(k_{A}, k_{A}) = k(-3n+1) \oplus k(-3n),
\]
with vector space basis $(b_{2} b_{3}^{n-1}, b_{3}^{n})$.  Theorem
\ref{thmA} gives us the following formulas in Ext:
\begin{align*}
m_{2}(b_{2} \otimes b_{3}) &= s_{5}, \\
m_{2}(b_{3} \otimes b_{2}) &= -s_{5}, \\
m_{2}(b_{3} \otimes b_{3}) &= s_{6}, \\
m_{3}(b_{2} \otimes b_{2} \otimes b_{2}) &= -s_{6}.
\end{align*}
All other instances of $m_{2}$ and $m_{3}$ on classes from $E^{1}$ are
zero, for degree reasons.
\end{example}

\begin{example}
Fix $p>2$ and let $A=k[x]/(x^{p})$ with $x$ in Adams degree $2d$ (so
that $A$ is graded commutative).  Its Ext algebra is
\[
\Ext_{A}^{*} (k_{A}, k_{A}) \cong \Lambda (y_{1}) \otimes k[y_{2}],
\]
with $y_{i}$ in $\Ext^{i}$, with $y_{1}$ in Adams degree $-2d$ and
$y_{2}$ in Adams degree $-2dp$.  Then Theorem \ref{thmA} tells us that
we may choose $y_{1}$ and $y_{2}$ so that $m_{p}(y_{1} \otimes \dotsb
\otimes y_{1}) = y_{2}$.  It is a standard Massey product computation
that the $p$-fold Massey product $\massey{y_{1}, \dotsc, y_{1}}$
equals a generator of $\Ext^{2}$ (with no indeterminacy, for degree
reasons), and Theorem \ref{thm-massey} tells us that with our choice
of $y_{1}$ and $y_{2}$, we have $\massey{y_{1}, \dotsc, y_{1}} =
\{(-1)^{(p+1)/2}y_{2}\}$.
\end{example}

\begin{example}
Let $k$ be a field of characteristic 2, and define the $k$-algebra $A$
by
\[
A = k\langle x_{1}, x_{2} \rangle / (x_{1}^{2}, \, x_{1} x_{2} x_{1} +
x_{2}^{2}),
\]
graded by putting $x_{i}$ in Adams degree $i$.  This is the sub-Hopf
algebra $A(1)$ of the mod 2 Steenrod algebra, and its cohomology can
be computed using spectral sequences -- see Wilkerson \cite[2.4]{Wi}
or Ravenel \cite[3.1.25]{Ra}, for instance.  In low degrees, it has
\[
\Ext^{n}_{A}(k_{A}, k_{A}) = \begin{cases}
k & n=0, \\
k(-1) \oplus k(-2), & n=1, \\
k(-2) \oplus k(-4), & n=2.
\end{cases}
\]
Following Wilkerson, we write $(h_{0}, h_{1})$ for the basis of
$\Ext^{1}$, with $h_{i}$ in Adams degree $-2^{i}$.  Then Theorem
\ref{thmA} tells us that $m_{2}(h_{0} \otimes h_{0}) \neq 0$ and
$m_{2}(h_{1} \otimes h_{1}) \neq 0$, so $(h_{0}^{2}, h_{1}^{2})$ is a
basis for $\Ext^{2}$.  Theorem \ref{thmA} also gives the formula
\begin{align*}
m_{3}(h_{0} \otimes h_{1} \otimes h_{0}) = h_{1}^{2}.
\end{align*}
This reflects the Massey product computation $\langle h_{0}, h_{1},
h_{0} \rangle = \{h_{1}^{2}\}$.
\end{example}

\section*{Acknowledgments} 

D.-M. Lu is supported by the Pao Yu-Kong and Pao Zhao-Long 
Scholarship and the NSFC (project 10571152). Q.-S. Wu is 
supported by the NSFC (project 10171016, key project 
10331030), in part by the grant of STCSM:03JC14013 in China, 
and by the Ky/Yu-Fen Fan Fund of the American Mathematical 
Society. J.J. Zhang is supported by grants from the National Science
Foundation (USA) and Leverhulme Research Interchange Grant 
F/00158/X (UK). A part of research was done when J.J. Zhang 
was visiting the Institute of Mathematics, Fudan University, 
China. J.J. Zhang thanks Fudan University for the warm 
hospitality and the financial support.

\end{document}